\documentclass{amsart}
\usepackage{hyperref}

\newtheorem{theorem}{Theorem}
\newtheorem{lemma}[theorem]{Lemma}

\newcommand{\R}{{\mathbb R}}
\newcommand{\Z}{{\mathbb Z}}

\newcommand{\be}{\begin{equation}}
\newcommand{\ee}{\end{equation}}

\newcommand{\ti}{\tilde}

\newcommand{\floor}[1]{\lfloor #1 \rfloor}

\newcommand{\I}{\mathrm{i}}

\newcommand{\Ran}{\mathrm{Ran}}

\newcommand{\eps}{\varepsilon}

\newcommand{\gam}{\gamma}

\numberwithin{equation}{section}

\begin{document}

\title[Bound States of Inverse Square Potentials]{Bound States of
Discrete Schr\"odinger Operators with Super-Critical Inverse Square Potentials}

\author[D. Damanik]{David Damanik}
\address{Mathematics 253--37\\ California Institute of Technology\\
Pasadena\\ CA 91125\\ U.S.A.}
\email{\href{mailto:damanik@caltech.edu}{damanik@caltech.edu}}
\urladdr{\href{http://www.math.caltech.edu/people/damanik.html}{http://www.math.caltech.edu/people/damanik.html}}

\author[G. Teschl]{Gerald Teschl}
\address{Fakult\"at f\"ur Mathematik\\
Nordbergstrasse 15\\ 1090 Wien\\ Austria\\ and International Erwin Schr\"odinger
Institute for Mathematical Physics, Boltzmanngasse 9\\ 1090 Wien\\ Austria}
\email{\href{mailto:Gerald.Teschl@univie.ac.at}{Gerald.Teschl@univie.ac.at}}
\urladdr{\href{http://www.mat.univie.ac.at/~gerald/}{http://www.mat.univie.ac.at/\~{}gerald/}}

\thanks{{\it Work supported by the National Science Foundation under Grant No.\ DMS-0500910
and the Austrian Science Fund (FWF) under Grant No.\ P17762}}


\begin{abstract}
We consider discrete one-dimensional Schr\"odinger operators whose potentials decay
asymptotically like an inverse square. In the super-critical case, where there are
infinitely many discrete eigenvalues, we compute precise asymptotics of the number of
eigenvalues below a given energy $E$ as this energy tends to the bottom of the essential
spectrum.
\end{abstract}

\maketitle

\section{Introduction}
This paper is concerned with discrete one-dimensional Schr\"odinger operators in
$\ell^2(\Z_+)$, where $\Z_+ = \{ 1,2,3,\ldots\ \}$. That is,
$$
H u(n) = - \Delta u(n) + V(n) u(n), \quad \Delta u(n) = u(n+1) - 2u(n) + u(n-1),
$$
where we impose a Dirichlet boundary condition, $u(0) = 0$.

If $V(n) \to 0$ as $n \to \infty$, then zero is the bottom of the
essential spectrum of $H$. We are interested in the discrete
spectrum of $H$ below zero. Thus, for $E \ge 0$, we define
$$
N_E(H) = \dim\Ran P_H((-\infty,-E]),
$$
where $P_H$ is the family of spectral projections associated with $H$ by the spectral
theorem. It is well known that $V(n) \sim - n^{-2+\eps}$ produces finitely many
eigenvalues if $\eps<0$ and infinitely many if $\eps>0$ and so inverse square decay is
critical for the existence of infinitely many discrete eigenvalues below the essential
spectrum. Furthermore, by the discrete analogue of Kneser's theorem, the discrete
spectrum below zero of the operator $H$ with potential $V(n) = - c n^{-2}$ is finite when
$c \le \frac14$ and infinite when $c > \frac14$ (see, e.g., \cite{dhs,lute,n}).

Our goal is to study the behavior of $N_E(H)$ as $E \downarrow 0$ in
the super-critical case $c>\frac14$. This question is also motivated
by recent results on a connection between singular spectrum embedded
in the essential spectrum and the discrete spectrum of a given
Schr\"odinger operator. See \cite{dhks,dk,dr,ks} and especially
\cite[Sect.~2]{dks}.

In the continuous case, Kirsch and Simon carried out an investigation of $N_E(H)$ for
super-critical inverse square potentials \cite{kisi} (see also \cite{schm} for
extensions). We prove the discrete analogue of their result. On the one hand, this case
is more relevant to the question raised in \cite[Sect.~2]{dks}. On the other hand, the
proof of Kirsch and Simon uses some arguments that do not carry over directly to the
discrete case: They scale the spatial variable and use exact solvability of the Euler
differential equation. Spatial scaling is not possible in the discrete case and, while
there exists a discrete Euler equation, it is not symmetric.

\begin{theorem}\label{main}
Suppose
$$
V(n) = - \frac{c}{n^2} + W(n), \qquad c> \frac{1}{4},
$$
where $W$ is a decaying sequence such that $N_0(-\Delta + \gam W)<\infty$ for all
$\gam\in\R$. Then
\begin{equation}\label{limit}
\lim_{E\downarrow 0} \frac{N_E(-\Delta+V)}{-\ln(E)} = \frac{1}{2\pi}
\sqrt{c-\tfrac{1}{4}}
\end{equation}
\end{theorem}
\noindent\textit{Remarks.} (i) We say that a sequence $W$ is decaying if $W(n) \to 0$
as $n \to \infty$.\\
(ii) The hypothesis on $W$ is satisfied, for example, if
$\sum_{n>0} n |W(n)|< \infty$. See \cite[Thm.~5.10]{tjac}.\\
(iii) The whole-line case can be reduced to the half-line case by Dirichlet decoupling.\\
(iv) For perturbations of the form $V(n) = \frac{c}{n^2} + W(n)$, an analogous result
holds near the top of the essential spectrum.

\section{Proof of Theorem~\ref{main}}

As a preparation we state the discrete analog of Proposition~5 from \cite{kisi}.
The proof is analogous.

\begin{lemma} \label{lemNE}
Let $V$, $W$ be decaying sequences. Then for every $E>0$ and $0<\eps<1$ we have
\begin{align*}
N_E\left(-\Delta + V + W\right) & \le N_E\left(-\Delta +\tfrac{1}{1-\eps} V\right) +
N_E\left(-\Delta + \tfrac{1}{\eps} W\right),\\
N_E\left(-\Delta + V + W\right) & \ge N_E\left(-\Delta + (1-\eps) V\right) - N_E
\left(-\Delta - \tfrac{1-\eps}{\eps} W\right).
\end{align*}
\end{lemma}
\noindent
Now we come to the proof of our main theorem. We start with $$V_c(n)= - \frac{c}{n^2}
$$ and replace it by $V_{E,c}$ which is just $V_c - E$ on $\{ n : V_c(n) \le -E\}$
and equal to $V$ otherwise. To investigate the asymptotics of $N_E(-\Delta + V_{E,c})$ we
split our domain into two parts by cutting at $\sqrt{\frac{c}{E}}$. For the first part,
we will compute the asymptotics of $N_E$ directly. The remaining part does not contribute
to $N_E$. Then we use Lemma~\ref{lemNE} to show that $N_E$ has the
same asymptotics for $V_{E,c}$ and $V= V_c + W$.

\begin{lemma}\label{lem3}
We have
$$
\lim_{E\downarrow 0} \frac{N_E(-\Delta+V_{E,c})}{-\ln(E)} = \frac{1}{2\pi}
\sqrt{c-\tfrac{1}{4}}.
$$
\end{lemma}

\begin{proof}
We first decompose $-\Delta+V_{E,c}$ into two parts by imposing an additional Dirichlet
boundary condition at $\floor{\sqrt{\frac{c}{E}}}$. Since this constitutes a rank-one
resolvent perturbation it will not affect the limit. By the choice of our cut point, the
part with $n>\floor{\sqrt{\frac{c}{E}}}$ does not contribute and by oscillation theory
(see e.g.\ \cite{tosc} or \cite[Ch.~4]{tjac}) it
suffices to count the number of sign flips of some solution of $(-\Delta + V_{E,c}) u =
-E u$ on $(1,\sqrt{\frac{c}{E}})$, that is, the number of sign flips of some solution of
$(-\Delta + V_c) u = 0$ on $(1,\sqrt{\frac{c}{E}})$.

Unfortunately, $(-\Delta + V_c) u = 0$ is not explicitly solvable, but
$$
\ti{u}_c(n)= \sqrt{n} \exp\left(\I \sqrt{c-\tfrac{1}{4}} \ln(n)\right)
$$
solves $(-\Delta + \ti{V}_c) \ti{u}=0$ with
the complex-valued potential
$$
\ti{V}_c(n) = \frac{\Delta \ti{u}_c(n)}{\ti{u}_c(n)} = -\frac{c}{n^2} + O(\frac{1}{n^3}).
$$
Moreover, it is straightforward to check
(cf.~\cite[Lemma~7.10]{tjac}, resp.~\cite{lute}) that $-\Delta u +
V_c u = 0$ has a solution $u_c$ which asymptotically looks like
$\ti{u}_c(n)$. Taking the real part of $u_c$, we see that the number
of sign flips behaves to leading order like $-\frac{1}{2\pi}
\sqrt{c-\frac{1}{4}} \ln(E)$.
\end{proof}
\noindent
Let us prove the upper bound in \eqref{limit}. By Lemma~\ref{lemNE},
\begin{align*}
N_E(-\Delta + V_c + W) &= N_E(-\Delta + (V_c - \chi_{(1-\eps)E,c}) +
N_E(\chi_{(1-\eps)E,c} + W))\\
&\le N_E(-\Delta + V_{E,c/(1-\eps)}) + N_E(-\Delta + \tfrac{1}{\eps} (\chi_{(1-\eps)E,c}
+ W)),
\end{align*}
where $\chi_{E,c}= E \chi_{(0,\sqrt{c/E})}$ and $\chi_\Omega$ is the characteristic function of the
set $\Omega$. Using
$$
N_E(-\Delta + \tfrac{1}{\eps} (\chi_{(1-\eps)E,c} + W)) \le N_0(-\Delta + \tfrac{1}{\eps}
W),
$$
the assumption on $W$, and Lemma~\ref{lem3}, we see that
$$
\limsup_{E\downarrow 0} \frac{N_E(-\Delta+V_c+W)}{-\ln(E)} \le \frac{1}{2\pi}
\sqrt{\tfrac{c}{1-\eps}-\tfrac{1}{4}}
$$
for every $0<\eps<1$, that is,
\begin{equation}\label{upper}
\limsup_{E\downarrow 0} \frac{N_E(-\Delta+V_c+W)}{-\ln(E)} \le \frac{1}{2\pi}
\sqrt{c-\tfrac{1}{4}}.
\end{equation}
It remains to show the lower bound in \eqref{limit}. By
Lemma~\ref{lemNE},
\begin{align*}
N_E(-\Delta + V_c + W) &= N_E(-\Delta + (V_c - \chi_{E/(1-\eps),c}) +
N_E(\chi_{E/(1-\eps),c} + W))\\
&\ge N_E(-\Delta + V_{E,(1-\eps) c}) - N_E(-\Delta - \tfrac{1-\eps}{\eps}
(\chi_{E/(1-\eps),c} + W))
\end{align*}
Observe that it suffices to show that the second
summand does not contribute to the limit. Invoking Lemma~\ref{lemNE} a second time we
have
$$
N_E(-\Delta - \tfrac{1-\eps}{\eps} (\chi_{E/(1-\eps),c} + W)) \le N_E(-\Delta -
\tfrac{1}{\eps} \chi_{E/(1-\eps),c}) + N_E(-\Delta - \tfrac{1-\eps}{\eps^2} W).
$$
The second term is bounded for fixed $\eps$ as $E\downarrow 0$ by assumption and it
remains to investigate the first one. As before we impose a Dirichlet boundary condition at
$\floor{\sqrt{\frac{c(1-\eps)}{E}}}$ and we need to count the sign flips of the solution
of $-\Delta u - \frac{E}{\eps(1-\eps)}u = -E u$ on $(0, \sqrt{\frac{c(1-\eps)}{E}})$.
Since this equation is explicitly solvable we obtain
$$
N_E(-\Delta - \tfrac{1}{\eps} \chi_{E/(1-\eps),c}) = \sqrt{c(1-\eps-\tfrac{1}{\eps})} + O(E).
$$
Hence
$$
\liminf_{E\downarrow 0} \frac{N_E(-\Delta+V_c+W)}{-\ln(E)} \ge
\frac{1}{2\pi} \sqrt{(1-\eps) c-\tfrac{1}{4}}
$$
for every $0 < \eps < 1$ and thus, \be \label{lower}
\liminf_{E\downarrow 0} \frac{N_E(-\Delta+V_c+W)}{-\ln(E)} \ge
\frac{1}{2\pi} \sqrt{c-\tfrac{1}{4}}. \ee Combining \eqref{upper}
and \eqref{lower}, we obtain the assertion of the theorem.
\hfill\qedsymbol

\bigskip

\noindent{\bf Acknowledgments.} G.T.\ gratefully acknowledges the
extraordinary hospitality of the Department of Mathematics at
Caltech, where this work was done.

\end{document}